\documentclass[reqno,12pt]{amsart2000}
\usepackage{amscd2000,amssymb}

\theoremstyle{plain}
\newtheorem{Thm}[subsection]{Theorem}
\newtheorem{Cor}[subsection]{Corollary}
\newtheorem{Lem}[subsection]{Lemma}
\newtheorem{Prop}[subsection]{Proposition}
\newtheorem{Conj}[subsection]{Conjecture}

\theoremstyle{definition}
\newtheorem{Def}[subsection]{Definition}

\theoremstyle{remark}

\newtheorem{Rem}[subsection]{Remark}

\newtheorem{conj}{Conjecture}

\numberwithin{equation}{section}
\renewcommand{\rm}{\normalshape}

\newif\ifShowLabels
\ShowLabelstrue
\newdimen\theight
\def\TeXref#1{%
        \leavevmode\vadjust{\setbox0=\hbox{{\tt
                \quad\quad  {\small \rm #1}}}%
        \theight=\ht0
        \advance\theight by \lineskip
        \kern -\theight \vbox to
        \theight{\rightline{\rlap{\box0}}%
        \vss}%
        }}%

\ShowLabelsfalse

\renewcommand{\sec}[2]{\section{#2}\label{S:#1}%
        \ifShowLabels \TeXref{{S:#1}} \fi}
\newcommand{\ssec}[2]{\subsection{#2}\label{SS:#1}%
        \ifShowLabels \TeXref{{SS:#1}} \fi}

\newcommand{\refs}[1]{Section ~\ref{S:#1}}
\newcommand{\refss}[1]{Section ~\ref{SS:#1}}

\newcommand{\reft}[1]{Theorem ~\ref{T:#1}}
\newcommand{\refl}[1]{Lemma ~\ref{L:#1}}
\newcommand{\refp}[1]{Proposition ~\ref{P:#1}}
\newcommand{\refc}[1]{Corollary ~\ref{C:#1}}

\newcommand{\refe}[1]{\eqref{E:#1}}
\newcommand{\refco}[1]{Conjecture ~\ref{Co:#1}}

\newenvironment{thm}[1]%
        { \begin{Thm} \label{T:#1}  \ifShowLabels \TeXref{T:#1} \fi }%
        { \end{Thm} }

\renewcommand{\th}[1]{\begin{thm}{#1} \sl }
\renewcommand{\eth}{\end{thm} }

\newenvironment{lemma}[1]%
        { \begin{Lem} \label{L:#1}  \ifShowLabels \TeXref{L:#1} \fi }%
        { \end{Lem} }
\newcommand{\lem}[1]{\begin{lemma}{#1} \sl}
\newcommand{\elem}{\end{lemma}}

\newenvironment{propos}[1]%
        { \begin{Prop} \label{P:#1}  \ifShowLabels \TeXref{P:#1} \fi }%
        { \end{Prop} }
\newcommand{\prop}[1]{\begin{propos}{#1}\sl }
\newcommand{\eprop}{\end{propos}}

\newenvironment{corol}[1]%
        { \begin{Cor} \label{C:#1}  \ifShowLabels \TeXref{C:#1} \fi }%
        { \end{Cor} }
\newcommand{\cor}[1]{\begin{corol}{#1} \sl }
\newcommand{\ecor}{\end{corol}}

\newenvironment{defeni}[1]%
        { \begin{Def} \label{D:#1}  \ifShowLabels \TeXref{D:#1} \fi }%
        { \end{Def} }
\newcommand{\defe}[1]{\begin{defeni}{#1} \sl }
\newcommand{\edefe}{\end{defeni}}

\newenvironment{remark}[1]%
        { \begin{Rem} \label{R:#1}  \ifShowLabels \TeXref{R:#1} \fi }%
        { \end{Rem} }
\newcommand{\rem}[1]{\begin{remark}{#1}}
\newcommand{\erem}{\end{remark}}

\newenvironment{conjec}[1]%
        { \begin{Conj} \label{Co:#1}  \ifShowLabels \TeXref{Co:#1} \fi }%
        { \end{Conj} }
\renewcommand{\conj}[1]{\begin{conjec}{#1} \sl }
\newcommand{\econj}{\end{conjec}}

\newcommand{\eq}[1]%
        { \ifShowLabels \TeXref{E:#1} \fi
           \begin{equation} \label{E:#1} }
\newcommand{\eeq}{ \end{equation} }

\newcommand{\prf}{ \begin{proof} }
\newcommand{\epr}{ \end{proof} }

\newcommand\alp{\alpha}         

\newcommand\gam{\gamma}         
\newcommand\del{\delta}         \newcommand\Del{\Delta}
\newcommand\eps{\varepsilon}

\newcommand\kap{\kappa}
\newcommand\lam{\lambda}                \newcommand\Lam{\Lambda}
\newcommand\sig{\sigma}         

\newcommand\ome{\omega}         

\newcommand\calC{{\mathcal{C}}}
\newcommand\calD{{\mathcal{D}}}

\newcommand\calF{{\mathcal{F}}}
\newcommand\calG{{\mathcal{G}}}
\newcommand\calH{{\mathcal{H}}}

\newcommand\calK{{\mathcal{K}}}

\newcommand\calM{{\mathcal{M}}}

\newcommand\calO{{\mathcal{O}}}
\newcommand\calP{{\mathcal{P}}}

\newcommand\calR{{\mathcal{R}}}
\newcommand\calS{{\mathcal{S}}}

\newcommand\calV{{\mathcal{V}}}
\newcommand\calW{{\mathcal{W}}}

            \newcommand\bfB{{\mathbf B}}
            
\newcommand\bfd{{\mathbf d}}

            \newcommand\bfG{{\mathbf G}}
            \newcommand\bfH{{\mathbf H}}

            \newcommand\bfM{{\mathbf M}}

\newcommand\bfp{{\mathbf p}}            \newcommand\bfP{{\mathbf P}}
\newcommand\bfq{{\mathbf q}}            \newcommand\bfQ{{\mathbf Q}}
            \newcommand\bfR{{\mathbf R}}
            
            \newcommand\bfT{{\mathbf T}}
\newcommand\bfu{{\mathbf u}}            \newcommand\bfU{{\mathbf U}}
            \newcommand\bfV{{\mathbf V}}
            
            \newcommand\bfX{{\mathbf X}}
            \newcommand\bfY{{\mathbf Y}}
            \newcommand\bfZ{{\mathbf Z}}

\newcommand\RR{\mathbb{R}}

\renewcommand\AA{\mathbb{A}}

\newcommand\FF{\mathbb{F}}
\newcommand\GG{\mathbb{G}}

\newcommand\ZZ{\mathbb{Z}}

\newcommand\CC{\mathbb{C}}

 \newcommand\grg{{\mathfrak{g}}}

 \newcommand\grm{{\mathfrak{m}}}

 \newcommand\grp{{\mathfrak{p}}}
 \newcommand\grq{{\mathfrak{q}}}
 \newcommand\grr{{\mathfrak{r}}}
 
 \newcommand\grt{{\mathfrak{t}}}
 \newcommand\gru{{\mathfrak{u}}}

\newcommand\sdp{\times \hskip -0.3em {\raise 0.3ex
\hbox{$\scriptscriptstyle |$}}} 

\newcommand\Hom{\operatorname {Hom}}

\newcommand\id{\operatorname{id}}
\newcommand\Id{\operatorname{Id}}

\newcommand\supp{\operatorname{supp}}

\newcommand\tr{\operatorname{tr}}

\newcommand\oP{{\overline{P}}}

\newcommand\oV{{\overline{V}}}

\newcommand\oX{{\overline{X}}}

\newcommand\opsi{{\overline{\psi}}}

\newcommand\hatF{{\widehat{F}}}

\newcommand\hatM{{\widehat{M}}}

\newcommand\hatT{{\widehat{T}}}

\newcommand\tilgam{{\widetilde{\gamma}}}

\newcommand{\obfp}{{\overline \bfp}}

\newcommand\x{\times}
\newcommand\ten{\otimes}

\newcommand{\ra}{\rangle}
\newcommand{\la}{\langle}

\renewcommand{\Id}{\text{Id}}

\newcommand\td{\bfT^{\vee}}

\newcommand\Gl{{\bf GL}}

\newcommand\Sp{{\bf Sp}}
\newcommand\Sl{{\bf SL}}

\newcommand\Av{\text{Av}}

\newcommand\ab{\text{ab}}
\newcommand\thsq{{\widetilde \bfH^2}}
\newcommand{\obfX}{{\overline \bfX}}
\newcommand{\obfP}{{\overline \bfP}}
\newcommand{\Gspin}{{\bf GSpin}}
\newcommand{\Spin}{{\bf Spin}}
\newcommand{\tbfP}{{\widetilde \bfP}}
\newcommand{\etj}{\eta_{P,Q,\psi}^{\text{junk}}}
\begin{document}
To Yuri Ivanovich Manin on the occasion of his 65th birthday

\

\

\

\title{Normalized intertwining operators and nilpotent elements in the Langlands dual
group}
\author{Alexander Braverman and David Kazhdan}
\address{Department of Mathematics, Harvard University, 1 Oxford st.
Cambridge MA, 02138}
\email{braval@math.harvard.edu, kazhdan@math.harvard.edu}
\begin{abstract}
Let $F$ be a local non-archimedian field and $\bfG$ be a split reductive
group over $F$ whose derived group is simply connected. Set $G=\bfG(F)$.
Let also
$\psi:F\to \CC^\x$ be a non-trivial additive unitary character of $F$.
For two parabolic subgroups $P$ and $Q$ in $G$ with the same Levi component
$M$ we construct an explicit unitary isomorphism
$\calF_{P,Q,\psi}:L^2(G/[P,P]){\widetilde\to}L^2(G/[Q,Q])$ commuting with the
natural actions of the group $G\x M/[M,M]$ on both sides. In some special
cases $\calF_{P,Q,\psi}$ is the standard Fourier transform. The crucial ingredient in
the definition is the action of the principal $sl_2$-subalgebra in the
Langlands dual Lie algebra $\grm^{\vee}$ on the nilpotent radical
$\gru_\grp^{\vee}$ of the Langlands dual parabolic.

For $M$ as above and using the operators $\calF_{P,Q,\psi}$
we define a {\it Schwartz space} $\calS(G,M)$. This space contains
the space $\calC_c(G/[P,P])$ of locally constant compactly supported
functions on $G/[P,P]$ for every $P$ for which $M$ is a Levi
component (but doesn't depend on $P$). We compute the space of spherical
vectors in $\calS(G,M)$ and study its global analogue.

Finally we apply the above results in order to give an alternative treatment
of automorphic L-functions associated with standard representations
of classical groups (thus reproving the results of \cite{GPSR} using the
same method as \cite{GJ}).

\end{abstract}
\maketitle

\sec{int}{Introduction}
\ssec{}{Notation}Let $F$ be a local non-archimedian field. Throughout this
paper we denote algebraic varieties over $F$ by boldface letters (e.g.
$\bfX$, $\bfG$ etc.) and their sets of $F$-points by the corresponding
ordinary letters ($G$, $X$ etc). If $\bfX$ is a smooth algebraic variety
over $F$ and $\ome$ is a top-degree differential form on $\bfX$, we denote
by $|\ome|$ the corresponding complex-valued measure on $X$ (cf. for example
\cite{weil}).

For $\bfX$ as above we denote by $\calC(X)$ the space of locally constant functions
on $X$. We also denote by $\calC_c(X)$ the space of locally constant functions with
compact support.

Let $\bfp:\bfX\to\bfY$ be a map of algebraic varieties. Then for a distribution $\eta$ on
$X$ we let $p_!(\eta)$ denote its direct image to $Y$ (it is well-defined if
the corresponding integral is convergent).
\ssec{}{Some results from \cite{BK99}}
The main part of this paper may be viewed as a continuation of \cite{BK99}.
In \cite{BK99} we considered the following situation. Let $\bfG$ be a split
reductive algebraic group over $F$ such that its derived group is simply
connected. Let also $\bfU\subset \bfG$ be the maximal unipotent subgroup
and set $\bfX=\bfG/\bfU$. Then $\bfX$ admits the natural action of
$\bfG\x\bfT$ where $\bfT$ is the (abstract) Cartan group of $\bfG$.

The variety $\bfX$ admits unique (up to multiplication by a constant)
$\bfG$-invariant top-degree differential form $dx$. Let $L^2(X)$ be the
space of $L^2$-functions on $X$ with respect to the measure $|dx|$.

In \cite{BK99} we considered (following Gelfand and Graev) the natural
unitary action of the Weyl group $W$ on the space $L^2(X)$. For every $w\in
W$ we denote by $\calF_w$ the corresponding operator.
The operators $\calF_w$ are given by certain explicit integral formulas
which generalize the standard Fourier transform. The main property of this
$W$-action is that it commutes with the natural $G$ action on $L^2(X)$ and
is compatible with the natural action of $W$ on $T$. Since the space
$L^2(X)$ is the direct integral of all representations induced from unitary
characters of a Borel subgroup $B$ it follows that one can think about
the operators $\calF_w$ as the "universal" {\it normalized intertwining
operators}.

In \cite{BK99} we also study the {\it Schwartz space} $\calS(X)$ of
functions on on $X$. By the definition, it is
the sum $\sum\limits_{w\in W}\calF_w(\calC_c(X))$
(the sum is taken in
$L^2(X)$). It is shown in \cite{BK99} that
$\calS(X)$ consists of locally constant functions. In \cite{BK99}
we also compute explicitly the space of spherical vectors in $\calS(X)$.
\ssec{}{Generalization to the parabolic case}In this paper we consider the
following generalization of the above results. Let $\bfM\subset \bfG$
be a Levi subgroup defined over $F$ and let $\bfM^\ab=\bfM/[\bfM,\bfM]$.
For every parabolic subgroup $\bfP\subset \bfG$
containing $\bfM$ we may consider the quotient
$\bfX_\bfP=\bfG/[\bfP,\bfP]$. Then $\bfX_\bfP$ is naturally a $\bfG\x
\bfM^\ab$-variety. Moreover, $\bfX_\bfP$ admits a $\bfG$-invariant
top-degree differential form which is unique up to
multiplication by a constant.
Hence it makes sense to consider $L^2(X_P)$.

The group $G$ acts naturally on $L^2(X_P)$ since it acts on $X_P$.
Let $\del_\bfP:\bfM^\ab\to\GG_m$ be the determinant of the action
of $\bfM$ on the Lie algebra of the unipotent radical $\bfU_\bfP$ of
$\bfP$ (by definition this action is the differential of $u\mapsto m^{-1}um$).
We define an action of $M^\ab$ on $L^2(X_P)$ by setting
\eq{normalization}
m(f)(x)=f(xm)|\del_\bfP(m)|^{\frac{1}{2}}
\end{equation}
Thus $f\mapsto m(f)$ is a unitary operator for each $m$.

Similarly, if $\eta$ is a distribution on $M^\ab$ then
we set
\eq{}
\eta(f)=\int\limits_{M^\ab} \eta(m)m(f).
\end{equation}
The following theorem is one of the main results of this paper (we are going to make it more precise
in the next subsection).
\th{existence}
Let $\psi:F\to \CC^\x$ be a non-trivial character. Then the following hold.
\begin{enumerate}
\item
Let $\bfP,\bfQ$ be any two parabolic subgroups of $\bfG$ which contain
$\bfM$ as their Levi factor. Then there exists canonical unitary isomorphism
\eq{}
\calF_{P,Q,\psi}:L^2(X_P)\to L^2(X_Q)
\end{equation}
commuting with the above $G\x M^\ab$-action.  Moreover, $\calF_{P,P,\psi}=\id$
\item
Let $\bfP,\bfQ,\bfR$ be three parabolic subgroups of $\bfG$ containing
$\bfM$ as their Levi factor. Then
$\calF_{Q,R,\psi}\circ\calF_{P,Q,\psi}=\calF_{P,R,\psi}$. Hence
$\calF_{P,Q,\psi}\circ\calF_{Q,P,\psi}=\id$.
\end{enumerate}
\eth
We shall omit the subscript $\psi$ in $\calF_{P,Q,\psi}$ when it does not lead
to a confusion.
\ssec{formula}{$\calF_{P,Q,\psi}$ and the action of the principal nilpotent in $\grm^{\vee}$}
Let us give the precise formula for the operator
$\calF_{P,Q,\psi}$. In this introduction we are going to ignore all convergence issues.
The rigorous treatment is given in \refs{torus} and \refs{proof}.

Recall that for $\bfP$ and $\bfQ$ as above there exists a non-normalized
intertwining operator $\calR_{P,Q}$ acting from $\calC_c(X_P)$ to
$\calC(X_Q)$ which commutes with the action of $G\x M^{\ab}$. This
operator is given by the following formula. Let $\bfZ_{P,Q}\subset
\bfX_\bfP\x\bfX_\bfQ$ be the image of $\bfG$ in $\bfX_\bfP\x\bfX_\bfQ$.
Then for every $f\in \calC_c(X_P)$
\eq{radon-parabolic}
\calR_{P,Q}(f)(y)=\int\limits_{(x,y)\in Z_{P,Q}}f(x)dx
\end{equation}

\noindent
{\bf Remark.} In fact, for \refe{radon-parabolic} to make sense we must choose $G$-invariant
measures on $G$ and $X_Q$ (both are well defined up to a constant). The choice of normalization
is explained in \refs{proof}.

We would like now to correct this operator in order to get a unitary
operator from $L^2(X_P)$ to $L^2(X_Q)$. More precisely we want
to construct a distribution $\eta_{P,Q,\psi}$ on $M^\ab$ such that
\eq{}
\calF_{P,Q,\psi}(f)=\eta_{P,Q,\psi}(\calR_{P,Q}(f)).
\end{equation}

Assume that $\bfT$ is a split torus over $F$. Let $\td$ denote the Langlands
dual torus. Set $\Lam_*(\bfT)=\Hom(\GG_m,\bfT)=\Hom(\td,\GG_m)$.
Let $L=\oplus L_i$ be a graded finite dimensional representation of
$\td$. Assuming that a certain technical condition on $L$ is satisfied
(cf. \refs{torus})  we can associate to $L$ a distribution $\eta_{L,\psi}$ on $T$ in the following way.
Choosing a homogeneous $\td$-eigen-basis for $L$ we may identify $L$ with
a collection $\lam_1,...,\lam_k$ of elements of $\Lam_*(\bfT)$ (with
multiplicities) with integers $n_1,...,n_k$ (the corresponding degrees) attached to them. Set
$s_i=\frac{n_i}{2}$. Then $\eta_{L,\psi}$ is the convolution of
distributions $(\lam_i)_!(\psi(t)|t|^{s_i}|dt|)$ (the technical condition mentioned above guarantees
the convergence of the above convolutions).

Let $\bfM^{\vee}$ be the Langlands dual group of $\bfM$. Then $(\bfM^\ab)^{\vee}=\bfZ(\bfM^\vee)$
(the center of $\bfM^\vee$). Let $\gru_\grp^\vee$ and $\gru_\bfq^\vee$ denote the nilpotent
radicals of the parabolic subalgebras of $\grg^{\vee}$ dual to $\bfP$ and $\bfQ$.
Let $\gru_{\grp,\grq}^\vee=\gru_\grp^{\vee}/\gru_\grp^\vee\cap\gru_\grq^{\vee}$.
Let $e,h,f$ be a principal $sl_2$-triple inside $\grm^\vee$. Consider $\gru_{P,Q}^e$. This is
naturally a representation of $\bfZ(\bfM^\vee)$ graded by the eigenvalues of $h$.
Set $\eta_{P,Q,\psi}=\eta_{(\gru_{\grp,\grq}^\vee)^e,\psi}$.

\th{formula}
\begin{enumerate}
\item
There exists a subspace $\calC_c^0(X_P)\subset\calC_c(X_P)$ which is dense in $L^2(X_P)$
such that for every
$f\in\calC_c(X_P)$ we have $\calR_{P,Q}(f)\in\calC_c(X_Q)$.
\item
For any $f\in \calC_c^0(X_P)$
define
$$
\calF_{P,Q,\psi}(f)=\eta_{P,Q,\psi}(R_{P,Q}(f)).
$$
Then $\calF_{P,Q,\psi}$ extends to a unitary operator from
$L^2(X_P)$ to $L^2(X_Q)$ which satisfies all the requirements of \reft{existence}.
\end{enumerate}
\eth
\ssec{intro-example}{Example}Let $V$ be a vector space over $F$ of dimension $n$ and
let $\oV$ denote the corresponding algebraic variety over $F$.
Consider $\bfG=\Sl(V)$. Let $\bfP\subset \bfG$ be the stabilizer of a line
$l$ in $V$ and let $Q$ be the stabilizer of a line in $V^*$ (which is in generic
position with respect to $l$). In this case $\bfX_\bfP=\bfV\backslash\{ 0\}$ and
$\bfX_\bfQ=\bfV^*\backslash\{ 0\}$. Hence we have
$$
L^2(X_P)=L^2(V)\quad\text{and}\quad L^2(X_Q)=L^2(V^*).
$$
In this case $\bfM^\ab=\GG_m$ and it
is easy to see that $\eta_{P,Q,\psi}(t)=\psi(t)|t|^{\frac{n-2}{2}}|dt|$
where $n=\dim\bfV$.
Also
$$
\bfZ_{\bfP,\bfQ}=\{ (v,v^*)\in\bfV\x\bfV^*|\ \la v,v^*\ra=1\}.
$$
It is easy to check that  $\calF_{P,Q,\psi}$ is just the Fourier transform (corresponding to
$\psi$) acting from $L^2(V)$ to $L^2(V^*)$.
\ssec{}{The space $\calS(G,M)$}For a Levy subgroup $M$ of $G$ let us denote by
$\calP(M)$ the set of all parabolic subgroups of $G$ containing $M$.
It follows from \reft{existence} that for any $P,Q\in\calP(M)$ we may identify the spaces
$L^2(X_P)$ and $L^2(X_Q)$. Thus we may regard all of them as one vector space which we shall
denote by $L^2(G,M)$. We define
$$
\calS(G,M)=\sum\limits_{P\in\calP(M)} \calC_c(X_P)\subset L^2(G,M).
$$
In the situation of
\refss{intro-example} one can prove that $\calS(G,M)$ is equal to the space
$\calC_c(V)$ (which is isomorphic to $\calC_c(V^*)$ by means of Fourier
transform). In the general case we don't have such a nice ``local'' description
of $\calS(G,M)$.

We study $\calS(G,M)$ in \refs{space} in some detail. In particular we compute its
subspace of spherical vectors (i.e. vectors invariant with respect to a standard
maximal compact subgroup).

\medskip
\noindent
{\it Remark.}
The formula for the spherical vectors in $\calS(G,M)$ is ``essentially equivalent''
to the formula for the intersection cohomology sheaf on Drinfeld's compactification
of the moduli space of $\bfP$-bundles on a smooth projective algebraic curve
(this intersection cohomology sheaf is studied in \cite{BFGM} (cf. also
\cite{FFKM} for the Borel case)). We don't have a geometric
explanation for this phenomenon (the main reason for this is that at the moment
we don't have an algebro-geometric analogue of the operators $\calF_{P,Q,\psi}$).
\ssec{}{Contents}This paper is organized as follows. In \refs{torus} we collect some auxiliary results
about distributions on a torus that will be used later. In \refs{borel} we recall some results
from \cite{GG}, \cite{K95} and \cite{BK99} and reformulate them in a little different
language. \refs{proof} is devoted to the proof of \reft{formula}. In \refs{space} we study the
{\it Schwartz space} $\calS(G,M)$ and compute its subspace of spherical vectors.
In \refs{global} we study the analogue of $\calS(G,M)$ when $F$ is replaced by a global field
$K$. We also formulate and prove certain analogue of the Poisson summation formula for the
operators $\calF_{P,Q,\psi}$ (in the situation of \refss{intro-example} it becomes the standard
Poisson summation formula for the Fourier transform).

Finally in \refs{L-function} we sketch how the above results may be applied in order to define and study
the (local and automorphic) L-functions associated with the standard representation of every classical group,
generalizing directly the method of \cite{GJ} where this is done for $\Gl(n)$. These $L$-functions
were studied in \cite{GPSR} by a different method and in \cite{PSR} by a method which is essentially equivalent
to ours.
However, the language of \cite{PSR} was more complicated since the operators $\calF_{P,Q,\psi}$ and
the above mentioned Poisson summation formula were not used there explicitly.
Thus one should think about \refs{L-function} as
a reformulation of \cite{PSR} using the results of this paper.
\ssec{}{Acknowledgments}We thank J.~Bernstein and S.~Rallis for very useful discussions on the
subject. We also thank Y.~Flicker for numerous comments about the text.
\sec{torus}{The distributions $\eta_{L,\psi}$}
\ssec{}{Good distributions on $F^\x$}
Let $\eta$ be a distribution on $F^\x$. We say that $\eta$ is good
if the following conditions hold:

1) For every open compact subgroup $K$ of $\calO^\x$ there exists a positive
real number $a=a(K)$ such that the integral
$$
\int_{t\in K} \eta(tx)
$$
vanishes for all $x$ such that $||x||>a$.

2) For every character $\chi:F^\x\to \CC^\x$ the Laurent power
series
$$
\sum\limits_{n=-\infty}^{\infty}z^n\int\limits_{||x||=q^{-n}} \eta(x)\chi(x)
$$
converges to a rational function $m_{\eta,\chi}(z)$. Note that the above series is indeed a Laurent
power series if we assume that 1 holds.

Let $\hatF^\x$ denote the group of all characters of $F^\x$. This set has a natural
structure of an algebraic variety over $\CC$ which is isomorphic to a disjoint
union of infinitely many copies of $\CC^\x$ indexed by
$\Hom(\calO^\x,\CC^\x)$.

Conditions 1 and 2 above imply that every good distribution $\eta$ on $F^\x$
defines a rational function $\calM(\eta)$ on $\hatF^\x$ (namely
for generic $\chi$ we have $\calM(\eta)(\chi)=m_{\eta,\chi}(1))$.
In this way $\calM$ becomes an isomorphism between the space of good
distributions on $F^\x$ and the space of rational functions on
$\hatF^\x$.

The basic example of a good distribution is the following. Let $\psi:F\to \CC^\x$
as before be a non-trivial additive character of $F$. Let also
$s$ be any complex number. Consider the distribution
$$
\eta_{\psi}^s=\psi(x)|x|^s|dx|.
$$
\lem{good}
$\eta_{\psi}^s$ is good.
\elem
The proof is left to the reader.

\medskip
\noindent
{\it Remark.} In this paper we are always going to normalize the measure $|dx|$ in such a
a way that the Fourier transform
$$
f(x)\mapsto \int f(x)\psi(xy)|dx|
$$
is a unitary operator (of course such normalization depends on the choice
of $\psi$).
\ssec{}{Good distributions on a torus}
We now generalize the above definitions to the case of an arbitrary split torus over
$F$. Let $\bfT$ be such a torus. We denote by $\Lam_*(\bfT)$ and $\Lam^*(\bfT)$ the corresponding
coweight and weight lattices. These lattices are naturally dual to each other. Note that $\bfT$ is canonically
defined over $\ZZ$ and $\bfT(\calO)$ is its maximal compact subgroup.

We have the natural valuation map
$v:T\to \Lam_*(\bfT)$ defined as follows. Let $t\in T$. Then for every
$\lam\in\Lam^*(\bfT)$ we define
$$
v(t)(\lam)=\text{the valuation of $\lam(t)\in F^\x$}.
$$
In this way $v$ defines an isomorphism between the quotient $T/\bfT(\calO)$ and $\Lam_*(\bfT)$.
For $\gam\in\Lam_*(\bfT)$ we set $T_\gam=v^{-1}(\gam)$.

Let $\grt_\RR=\Lam_*(\bfT)\ten \RR$ and let $\calK\in\grt_\RR$
be a closed cone satisfying the following conditions:

a) $\calK$ is generated by finitely many elements of $\Lam_*(\bfT)$.

b) The interior of $\calK$ is open in $\grt_\RR$.

c) $\calK$ contains no straight lines.

Let also $\hatT=\Hom(T,\CC^\x)$. Then $\hatT$ has a natural
structure of an algebraic variety over $\CC$ which is isomorphic
to a disjoint union of infinitely many copies $(\CC^*)^{\dim \bfT}$
indexed by $\Hom(T(\calO),\CC^\x)$

Let $\eta$ be a distribution on $T$. We say that $\eta$ is $\calK$-good
if the the following conditions hold:

1) For every open compact subgroup $K$ of $T$ there exists $\gam_0\in\Lam_*(\bfT)$ such that the integral
$$
\int_{t\in K}\eta(tx)
$$
vanishes for every $x\in T$ such that $v(x)\not\in \gam_0+\calK$. We denote
by $\Av_K(\eta)$ the $K$-invariant distribution on
$T$ given by the above integral (where $x$ is considered as a variable). We also denote
by $\Av_K(\eta)_\gam$ its restriction to $T_\gam$.

2) $\Av_K(\eta)$ has polynomial growth, i.e. there exists a polynomial function $p$
$\grt_\RR$  and a functional $\lam\in\grt_\RR^*$ such that for every $\gam\in\Lam_*(\bfT)$ we have
$$
||\Av_K(\eta)_\gam||_{L^1}\leq |p(\gam)|q^{\lam(\gam)},
$$
where $||\Av_K(\eta)||_{L^1}$ denotes the $L^1$-norm of $\Av_K(\eta)_\gam$.

Condition 2
implies that there exists an open subset of every connected component of $\hatT$ such that
for every $\chi$ in this subset the integral
$$
\int_T \eta(t)\chi(t)=\sum_{\gam\in \Lam_*(\bfT)}\int_{T_\gam}\eta(t)\chi(t)
$$
is absolutely convergent. Thus we can impose the following (last) condition
on $\eta$:

3) $\eta$ defines a rational function on
$\hatT$, i.e. there exists a rational function $\calM(\eta)$ on
$\hatT$ such that for
$$
\calM(\eta)(\chi)=\int_T \eta(t)\chi(t)
$$
if the latter integral is absolutely convergent.

We denote the space of $\calK$-good distributions by $\calD_\calK(T)$.

Given $\eta_1,\eta_2\in\calD_\sig(T)$ the convolution $\eta_1*\eta_2$ makes sense.
Indeed, for every open compact subgroup $K$ of $T$ the convolution
$\Av_K(\eta_1)*\Av_K(\eta_2)$ makes sense because it is defined by
a proper integral.
Let now $\phi$ be a test function, i.e. a compactly supported function
on $T$ which is invariant under some maximal compact subgroup $K\subset T$.
Then
$$
(\eta_1*\eta_2)(\phi)=\Av_K(\eta_1)*\Av_K(\eta_2)(\phi)
$$
 It is easy to see that
$$
\calM(\eta_1*\eta_2)=\calM(\eta_1)\cdot\calM(\eta_2).
$$
It is clear from this formula that $\eta_1*\eta_2\in\calD_\calK(T)$.
\ssec{maindist}{Example}Let $\lam_1,...,\lam_k$ be a collection of non-zero elements of $\calK$.
Let also $s_1,...,s_k$ be some complex
numbers. In this case we define the distribution
$$
\eta^{s_1,...,s_k}_{\lam_1,...,\lam_k,\psi}=(\lam_1)_!(\eta_{\psi,s_1})*...*(\lam_k)_!(\eta_{\psi,s_k}).
$$
If $\bfT=\GG_m$ and all $\lam_i$ are equal to the standard character of $\GG_m$ then we shall
just write $\eta_\psi^{s_1,...,s_k}$ for the above distribution.

Let $L=\oplus L_i$ be a graded representation of $\td$.
Choosing a homogeneous $\td$-eigen-basis we may identify $L$
with a collection $\lam_1,...,\lam_k$ of elements of $\Lam_*(\bfT)$ with
integers $n_1,...,n_k$ attached to them. Set $s_j=\frac{n_j}{2}$.

We say that $L$ is $\calK$-good if $\lam_i\in\calK$ for every $i$.
In this case we define the distribution $\eta_{L,\psi}$ by setting
\eq{}
\eta_{L,\psi}=\eta^{s_1,...,s_k}_{\lam_1,...,\lam_k,\psi}
\end{equation}
It follows from \refl{good} and from the fact that $L$ is $\calK$-good that
$\eta_{L,\psi}$ is well-defined and belongs to $\calD_\calK(T)$.
\ssec{}{Unitarity properties}Let $\bfX$ be an algebraic variety over $F$ endowed with
a free $\bfT$-action. Let $dx$ be top-degree differential form on $\bfX$.
Assume that there exists a character $\del:\bfT\to\GG_m$ such that
$d(t^{-1}\cdot x)=\del(t)dx$. In this case we define an action of $T$ on functions
on $X$ by setting
$$
t(f)(x)=|\del(t)|^{1/2}f(t^{-1}x).
$$
Thus every $t\in T$ acts on the space $L^2(X,|dx|)$ as a unitary operator.

As before for a distribution $\eta$ on $T$ we write
$$
\eta(f)=\int_T t(f)\eta(t).
$$

In what follows we say that some statement holds for generic
$f\in \calC_c(X)$ if it there exists a subspace
$\calC_c^0(X)$ which is dense in $L^2(X)$ such that
the above statement holds for any $f\in\calC_c^0(X)$.

We say that a distribution $\eta$ is {\it unitary} if for generic
$f\in\calC_c(X)$ we have
$$
||\eta(f)||=||f||.
$$

\lem{gm}
Let $\bfT=\GG_m$. Then $\eta_\psi^{s}$ is unitary  if and
only if $s=-\frac{1}{2}$. Also, $\eta_\psi^{s_1,s_2}$ is unitary on $L^2(X)$
if and only if $s_1+s_2=-1$.
\elem
\prf
Because of our normalizations it is enough to assume that $\bfX=\GG_m$ with the multiplicative
measure $d^*x=\frac{|dx|}{|x|}$ on $X$. Thus we have
$$
\eta_\psi^s(x)=\psi(x)|x|^s |dx|=\psi(x)|x|^{s+1}d^*x.
$$
Let $A_s$ denote the operator of convolution with $\eta_\psi^s$.
Define also an operator $B:\calC_c(F^\x)\to \calC(F^\x)$ by
$$
B_s(f)(y)=\int\limits_{F^\x} f(x)\opsi(xy^{-1})\frac{|x|^{s+1}}{|y|^{s+1}}d^*x.
$$
It is easy to see that $B_s$ is the Hermitian conjugate of $A_s$. More precisely we claim
that there exists a subspace $\calC_c^0(F^\x)$ such that

a) $\calC_c^0(F^\x)$ which is dense in $L^2(F^\x)$

b) For every $f\in\calC_c^0(F^\x)$ both
$A_s(f)$ and $B_s(f)$ lie in $\calC_c(F^\x)$.

c) For any $f,g\in\calC_c^0(F^\x)$ we have
$$
\la A_s(f),g\ra=\la f,B_s(g)\ra.
$$

To prove \refl{gm} it is enough to show that that
$B_{s_2}\circ A_{s_1}=\Id$ if $s_1+s_2=-1$ (the ``only if''
statement thus follows automatically).

However,

\begin{align}
B_{s_2}\circ A_{s_1}(f)(z)=\int
f(x)\psi(x^{-1}y)|x^{-1}y|^{s_1+1}\opsi(yz^{-1})|yz^{-1}|^{s_2+1}d^*x d^* y=\\
\int f(x)\psi(y(x^{-1}-z^{-1}))|y|^{s_1+s_2+1}|dy||x|^{-s_1-1}|z|^{-s_2-1}d^*x
\end{align}
The latter integral is clearly equal $f(z)$ if $s_1+s_2=-1$.
\epr
\cor{unitarity}
Let $\bfX$ be as above. Let $\lam\in\Lam_*(\bfT)$ be any non-zero element.
Then $\eta_{\lam,\psi}^s$ is unitary if and
only if $s=-\frac{1}{2}$. Also $\eta_{\lam,\lam,\psi}^{s_1,s_2}$ is unitary
if and only if $s_1+s_2=-1$.
\ecor


\cor{main-unitarity}
Let $L'$ be a finite-dimensional $sl(2)\x \td$-module and and $L$ be its
quotient by the space of highest weight vectors. Consider the graiding
on $L$ by the eigenvalues of $h$. Then $\eta_{L,\psi}$
induces a unitary operator on the space $L^2(X,|dx|)$.
\ecor
\prf It is easy to see that for every $\lam\in\Lam_*(\bfT)$ the multiplicity of  $n$ as an eigenvalue
of $h$ in $L_\lam$ is equal to the multiplicity of $2-n$. Thus the proof follows
from \refc{unitarity}.
\epr
\ssec{radon}{Example: Fourier transform vs. Radon transform}
Let $\bfp:\bfX\to \bfY$ be a rank $n$ vector bundle over a smooth variety
$\bfY$ defined over $F$ and let $\bfp^{\vee}:\bfX^{\vee}\to \bfY$ be the dual vector
bundle. Let $dy$ be a volume form on $\bfY$. Let $\ome$ be a non-vanishing
section of $\det \bfX^{\vee}$ on $\bfY$ and let $\ome^{\vee}$ be the corresponding
section of $\det X$. Then $\ome\ten dy$ makes sense as a volume
form on $\bfX$ which we denote by $dx$. Similarly we define $dx^{\vee}$.
We have the Fourier transform $\calF_X:L^2(X,|dx|)\to L^2(X^\vee, |dx^\vee|)$.

Let $\bfX_0$ (resp. $\bfX^{\vee}_0$) be the complement to the zero section
in $\bfX$ (resp. in $\bfX^\vee$).
Note that $L^2(X,|dx|)=L^2(X_0,|dx|)$. We let the group $\GG_m$ act freely on
$\bfX_0$ and $\bfX^\vee_0$ by setting
$$
t\cdot x=t^{-1}x\quad t\cdot x^\vee=tx^{\vee}
$$
(where in the right hand side of the above equalities we mean usual action
of scalars on the fibers of a vector bundle).
The reasons for such normalization will (hopefully) become clear in the next
section.

Let
$$
\bfZ_X=\{ (x,x^{\vee})\in \bfX\underset{\bfY}\x\bfX^\vee|\quad \la x,x^\vee\ra=1\}.
$$
and let $\pi:\bfZ_\bfX\to\bfX,\pi^\vee:\bfZ_\bfX\to\bfX^{\vee}$
be the natural projections. There is a natural fiberwise volume
form along the fibers of either $\pi$ or $\pi^\vee$. Hence the operations $\pi_!$ and
$\pi^\vee_!$ make sense on compactly supported functions.

For $f\in\calC_c(X_0)$ let
$$
\calR_X(f)=\pi^\vee_! \pi^*(f).
$$
It is easy to see that the intersection of the support of $\calR_X(f)$
with every $F^\x$-orbit is compact. Hence the convolution of $\calR_X(f)$ with any
distribution on $F^\x$ is well-defined.

Let $\eta(t)=\psi(t)|t|^{\frac{n-2}{2}}|dt|$.
\lem{fourier-radon}
$$
\calF_X(f)=\eta(\calR_X(f)).
$$
\elem

\sec{borel}{Digression on the Borel case}
In this section we review the case when $\bfM=\bfT$ is a maximal split torus
in $\bfG$, i.e. all the parabolic subgroups in question are Borel subgroups.
In this case the formula for $\calF_{P,Q}$ is essentially due to Gelfand and
Graev.

Let $\bfB,\bfB'$ be two Borel subgroups of $\bfG$. As is well-known to such
a pair one can canonically associate an element $w$ of the Weyl group.
It is clear that the collection of operators $\calF_{B,B'}$ satisfying
the conditions of \reft{existence} is uniquely determined by
those for which $w$ is a simple reflection. Below we give a formula
for $\calF_{B,B'}$ in that case.

\ssec{}{The case of a simple reflection}The construction explained below is a reformulation of the construction
of \cite{GG} (cf. also \cite{K95} and \cite{BK99}).

Let $\alp$ be a simple root of $\bfG$ and assume that $B$ and $B'$ are
in position $w=s_\alp$. Let $\bfP$ be the minimal parabolic subgroup containing
both $\bfB$ and $\bfB'$. Let $\bfp:\bfX_\bfB\to\bfX_\bfP$ (resp. $\bfp':\bfX_{\bfB'}\to\bfX_\bfP$) be the
natural projections. Let also $\obfp:\obfX_\bfB\to\bfX_\bfP$ (resp.
$\obfp':\obfX_{\bfB'}\to\bfX_\bfP$)
be their affine completions (i.e. $\obfp$ is the affine morphism
corresponding to the sheaf of algebras $\bfp_*\calO_{\bfX_\bfB}$ on
$\bfX_\bfP$). Then $\obfp$ and $\obfp'$ are mutually dual vector bundles
over $\bfX_\bfP$. Let us explain how the natural pairing
$\kap:\bfX_\bfB\underset{\bfX_\bfP}\x\bfX_\bfB'\to \AA^1$ looks like.
Let $\alp^\vee:\GG_m\to \bfT$ be the simple coroot corresponding to $\alp$.
Then $\kap$ is uniquely characterized by the following two requirements:

\noindent
1) For every $g\in G$ we have
\eq{kap1}
\kap(g\text{mod} \bfU,g\text{mod}\bfU')=1.
\end{equation}
2) For every $(x,y)\in \bfX_\bfB\underset{\bfX_\bfP}\x\bfX_\bfB'$ and every $t\in\bfT$
one has
\eq{kap2}
\kap(\alp^\vee(t)x,y)=\kap(x,\alp^\vee(t^{-1}y))=t\cdot\kap(x,y).
\end{equation}

In what follows we choose $G$-invariant measures on $G$ and $X_B$ (for all Borel subgroups
$\bfB$ defined over $F$) in such a way that the measure of the image of $\bfG(\calO)$ is
equal to 1. This measure extends naturally to $\oX_B$.

Note that $L^2(X_B)=L^2(\oX_B)$.
We define $\calF_{B,B',\psi}$ to be the Fourier transform in the fibers of
the bundle $\obfp$. It makes sense as a unitary operator acting from $L^2(X_B)$ to
$L^2(X_{B'})$.
\ssec{}{The general case}Let $\bfB,\bfB'$ be any two Borel subgroups of $\bfG$.
Then there exists a sequence $(\bfB_0=\bfB,\bfB_1,...,\bfB_n=\bfB')$ of Borel subgroups
defined over $F$
such that $\bfB_i$ and $\bfB_{i+1}$ are in position $s_{\alp_i}$ where
$\alp_i$ is a simple root of $\bfG$ (for every $i$).
We define
$$
\calF_{B,B',\psi}=\calF_{B_{n-1},B_n,\psi}\circ ...\circ \calF_{B_0,B_1,\psi}.
$$
It is shown in \cite{K95} that $\calF_{B,B',\psi}$ does not depend on the choice
of the sequence $(\bfB_0,\bfB_1,...,\bfB_n)$. It is clear that the operators
$\calF_{B,B',\psi}$ are unitary and satisfy the requirements of \reft{existence}.
\ssec{}{The distribution $\eta_{B,B',\psi}$}
Let $\bfB,\bfB'$ be as above and assume that they are in position $w\in W$.
We choose a maximal torus $\bfT\subset \bfB\cap\bfB'$.

Let $\Pi_w$ denote the set of all positive with respect to $\bfB$ coroots of
$\bfT$ which are made negative by $w$. Let also $\calK\subset\grt_\RR$ be
the cone of positive coroots.

Thus for every $\alp\in\Pi_V$ the distribution $\eta_{\alp,\psi}^0$ is $\calK$-good.
Therefore, their convolution makes sense. We define
$\eta_{B,B',\psi}$ to be equal to the convolution of $\eta_{\alp,\psi}^0$ for all
$\alp\in\Pi_w$. The reader will readily check that this definition of
$\eta_{B,B',\psi}$ coincides with the one given in \refss{formula}.
\ssec{}{Proof of \reft{formula} in the Borel case}Let $\bfB$, $\bfB'$ as before be two
Borel subgroups defined over $F$. Thus we may consider the operator $\calR_{B,B'}$
defined by \refe{radon-parabolic}. The operator $\calR_{B,B'}$ is well defined as an operator
from $\calC_c(X_B)$ to $\calC(X_{B'})$. Moreover, it is easy to see that
for every $f\in \calC_c(X_B)$ the intersection of $\text{supp}\calR_{B,B'}(f)$ with
any orbit of $T$ is compact. Thus $\eta(\calR_{B,B'}(f))$ is well-defined for any
distribution $\eta$ on $T$. We claim that for every $f\in \calC_c(X_B)$ we have
\eq{dostalo}
\calF_{B,B',\psi}(f)=\eta_{B,B',\psi}(\calR_{B,B'}(f)).
\end{equation}
Let $B,B',B''$ be three Borel subgroups. Let $w_1$ be the relative position of $B$ and $B'$
and let $w_2$ be the relative position of $B'$ and $B''$. Assume that
$l(w_2w_1)=l(w_2)+l(w_1)$ (here $l$ denotes the length function on $W$).
It is easy to see that in this case we have
$$
\calR_{B',B''}\circ\calR_{B,B'}=\calR_{B,B''}\quad\text{and}\quad
\eta_{B',B'',\psi}\star \eta_{B,B',\psi}=\eta_{B,B'',\psi}
$$
Hence it is enough to check \refe{dostalo} when the relative position of $B$ and $B'$
is a simple reflection. In this case it follows from \refl{fourier-radon}.
\sec{proof}{Proof of \reft{existence} and \reft{formula}}
\ssec{}{}Let $M$ be a Levi subgroup of $G$ and let $P,Q\in \calP(M)$. Choose a Borel subgroup
$B$ contained in $P$ and let $B'$ another Borel subgroup such that:

1) $B'$ is contained in $Q$.

2) The relative position $w$ of $B$ and $B'$ has minimal length subject to the first requirement.

Let $\pi:X_B\to X_P$, $\pi':X_{B'}\to X_Q$ denote the natural projections. Since all the spaces
in question are endowed with natural measures the operations $\pi_!$ and $\pi'_!$ may be applied
to functions.
\prop{junk}
Recall that $\gru_{\grp,\grq}^\vee=\gru_\grp^\vee/\gru_\grp^\vee\cap\gru_\grq^\vee$; this space
is graded by the eignevalues of $h$ where $(e,h,f)\in \grm^\vee$ is a principal $sl_2$-triple.
\begin{enumerate}
\item For generic $\phi\in\calC_c(X_P)$ we have
$\calR_{P,Q}(\phi)\in\calC_c(X_Q)$ and
$\calR_{Q,P}\circ\calR_{P,Q}(\phi)\in\calC_c(X_P)$.
\item
There exists unique unitary operator
$\calG_{P,Q,\psi}:L^2(X_P)\to L^2(X_Q)$ such that for every
$f\in \calC_c(X_B)$ we have
\eq{dostalo2}
\calG_{P,Q,\psi}(\pi_!(f))=\pi'_!(\calF_{B,B',\psi}(f)).
\end{equation}
\item
We have
\eq{dostalo3}
\calG_{P,Q,\psi}(\phi)=\eta_{\gru_{\grp,\grq}^\vee,\psi}(\calR_{P,Q}(\phi))
\end{equation}
for generic $\phi\in\calC_c(X_P)$.
\end{enumerate}
\eprop
\prf
Let us prove 1. For a character $\chi:M^\ab\to\CC^\x$ we denote by
$\calC_c(X_P)_\chi$ the corresponding space of $(M,\chi)$-coinvariants
in $\calC_c(X_P)$ (this space is dual to the corresponding
induced representation). It is well-known that for generic (i.e. lying
in a Zariski dense subset) $\chi$ the operators $\calR_{P,Q}$ and $\calR_{Q,P}$
give rise to isomorphisms $\calC_c(X_P)_\chi\simeq\calC_c(X_Q)_\chi$.
Let now $D\subset \hatM^\ab$ denote the complement to the above
Zariski dense subset. Then it is clear that 1 holds
for every $\phi$ whose image in $\calC_c(X_P)_\chi$ is
equal to $0$ for every $\chi\in D$.

Let us now prove 2 and 3.
Let $\phi\in\calC_c(X_P)$. Then there exists a function $f\in \calC_c(X_B)$
such that $\pi_!(f)=\phi$. Therefore $\calG_{P,Q,\psi}(\phi)$ is uniquely
defined by \refe{dostalo2}.
On the other hand it is easy to see that for every $\phi\in\calC_c^0(X_P)$ and $f$
as above  the right hand side of
\refe{dostalo2} is equal to \refe{dostalo3}. This shows that $\calG_{P,Q,\psi}(\phi)$
is well-defined as an operator $\calC_c^0(X_P)\to \calC(X_Q)$.

Let $\phi\in\calC_c(X_P), \phi'\in\calC_c(X_Q)$. Then it follows from
\refe{dostalo} that
$$
\la \calG_{P,Q,\psi}(\phi),\phi'\ra=\la \phi,\calG_{Q,P,\opsi}(\phi')\ra.
$$
 On the other hand,
since $\calF_{B,B',\psi}^{-1}=\calF_{B',B,\opsi}$ it follows that
that $\calG_{P,Q,\psi}^{-1}=\calG_{Q,P,\opsi}$. Thus the inverse
of $\calG_{P,Q,\psi}$ is equal to its hermitian conjugate which means that
$\calG_{P,Q,\psi}$ is unitary.
\epr
\ssec{}{Proof of \reft{formula}}We have to prove that the operator $\calF_{P,Q,\psi}$ defined
as in \reft{formula} is unitary.

Let $(\gru_{\grp,\grq}^\vee)^{\text{junk}}=\gru_{\grp,\grq}^\vee/(\gru_{\grp,\grq}^\vee)^e$
(with grading induced by the eigenvalues of $h$). Let
$\etj=\eta_{(\gru_{\grp,\grq}^\vee)^{\text{junk}},\psi}$. Thus
$$
\eta_{\gru_{\grp,\grq}^\vee,\psi}=\eta_{P,Q,\psi}\star \etj.
$$
Hence it follows from \refp{junk} that for every $\phi\in \calC_c(X_P)$ we have
$$
\calG_{P,Q,\psi}=\etj(\calF_{P,Q,\psi}(\phi)).
$$
Since $\calG_{P,Q,\psi}$ is a unitary operator it follows that to prove the unitarity
of $\calF_{P,Q,\psi}$ it is enough to prove the unitarity of the operator of convolution
with $\etj$. This, however, follows immediately from \refc{main-unitarity}.

We have to prove now the assertion of Theorem 1.4(2). However, it is enough
to do it in the following 2 cases:

1) $R=P$.

2) $\dim \gru_{\grp,\grr}^\vee=\dim\gru_{\grp,\grq}+\dim\gru_{\grq,\grr}$.

In case 1 the statement follows immediately from the unitarity of
$\calF_{P,Q,\psi}$ and in case 2 this is obvious from the definitions.
\ssec{glnexample1}{Example}Consider the case when $\bfG=\Sl(n)$ and
$\bfM=\Gl(n-1)$ embedded into $\bfG$ in the standard way:
$$
(a_{ij})_{i,j=1}^n\mapsto
\begin{pmatrix}
a_{11}\quad...\quad a_{1,n-1}& 0\\
...&\\
a_{n-1,1}\quad ...\quad a_{n-1,n-1}& 0\\
0 \quad ...\quad 0& \det(a_{ij})^{-1}
\end{pmatrix}
$$
There are two parabolic subgroups $\bfP$ and $\obfP$ containing $\bfM$.
Namely, we set
$$
\bfP=\text{stab}(\text{span}
\begin{pmatrix}
1\\
0\\
...\\
0
\end{pmatrix})
$$
and $\obfP$ is the corresponding opposite parabolic.

Let $V$ denote the defining representation of $G$.
Then $X_P=V\backslash\{0\}$ and $X_{\oP}=V^*\backslash\{0\}$.
Hence we have $L^2(X_P)=L^2(V)$ and $L^2(V^*)$ (where the measure on $V$ and
$V^*$ is a Haar measure with respect to addition).
It follows from \refl{fourier-radon} that $\calF_{P,\oP,\psi}$ is equal to
the Fourier transform $\calF_{V,\psi}$.


\sec{space}{The space $\calS(G,M)$}
In this section we assume that the character $\psi$ is trivial on $\calO$ and that it is non-trivial
on $\pi^{-1}\calO$.
\ssec{}{}Let $M\subset G$ be a Levi subgroup. It follows from the result of the previous section
that for every two parabolic subgroups $P$ and $Q$ for which $M$ is the Levi factor
the spaces $L^2(X_P)$ and $L^2(X_Q)$ are canonically isomorphic.
Hence we may regard it as one space which we shall denote by $L^2(G,M)$.

For every $P$ as above we denote by $C_c(X_P)$ the space of compactly supported locally constant
functions on $X_P$. We have the natural embedding $C_c(X_P)\subset L^2(G,M)$.
We define
\eq{}
\calS(G,M)=\sum \calC_c(X_P)\subset L^2(G,M)
\end{equation}
(the sum is being taken over all parabolic subgroups in which $M$ is a Levi factor).
Clearly, $\calS_(G,M)$ is a representation of $G\x M^\ab$. In the case $M=T$
this space has been studied in \cite{BK99}.
\ssec{glnexample}{Example}Assume that we are in the situation
of \refss{glnexample1}.
\lem{schwartz-gln}
In this case we have
$$
\calS(G,M)=\calC_c(V).
$$
\elem
\prf
Since $\calC_c(V)$ is invariant under Fourier transform, it follows that
$\calS(G,M)$ is a subspace of $\calC_c(G,M)$. Since
$\calC_c(X_P)=\calC_c(V\backslash\{0\}$ has codimension 1 in $\calC_c(V)$,
in order to prove the opposite inclusion it is enough to find $f\in\calC_c(V^*\backslash\{0\}$
whose Fourier transform does not vanish at $0$; we may take
any $f$ whose integral over $V^*$ does not vanish
Let $K=\bfG(\calO)$ be the standard maximal compact subgroup of $G$.
\epr

In the case of arbitrary $G$ and $M$ we don't know any "nice" description of
$\calS(G,M)$. Let us, however, discuss some simple properties of its
elements.

For every $f\in\calS(G,M)$ we denote by $f_P$ the corresponding
function on $X_P$.
\lem{locallyconstant}Let $f\in \calS(G,M)$. Then $f_P$ is a locally constant function
on $X_P$ for every $P\in\calP(M)$.
\elem
\prf
Clearly it is enough to show the following:
let $P,Q\in \calP(M)$ and let $h\in \calC_c(X_P)$. Then
$\calF_{P,Q,\psi}(h)$ is locally constant on $X_Q$.

To prove this let us note that $h$ is fixed by some compact subgroup $C$ of $G$. Since
$\calF_{P,Q}$ commutes with the action of $G$ the same is true for $\calF_{P,Q}(h)$.
Since $G$ acts transitively on $X_Q$ this implies that $h$ is locally constant.
\epr

Let $\Lam_*(\bfM)$ denote the lattice of cocharacters
of $\bfM^{\ab}$. Let also $\Lam_*=\Lam_*(\bfT)$ be the coroot lattice of $\bfG$.
 We have the natural restriction map $\Lam_*\to\Lam_*(\bfM)$.

Fix $\gam\in\Lam_*(\bfM)$. Let $\tilgam$ be any lift of $\gam$ to an element of
$\Lam_*$. It is easy to see that the $K$-orbit of $\tilgam(\pi)\text{mod}[P,P]$ in $X_P$ depends only
on $\gam$. We denote this orbit by $X_P^\gam$.
The following lemma is well-known:
\lem{}
The assignment $\gam\mapsto X_P^\gam$ is a one-to-one correspondence between
$\Lam_*(\bfM)$ and the set of $K$-orbits on $X_P$.
\elem

Let $\Lam_*^+$ denote the set of all
linear combinations of positive coroots of $\bfG$ with non-negative coefficients. We say that
$\gam\in\Lam_*(\bfM)$ is positive if it is equal to the image of some element of $\Lam_*^+$.
We denote by $\Lam_*(\bfM)^+$ the set of positive elements in $\Lam_*(\bfM)$. We say that
a function $f$ on $X_P$ has {\it bounded support} if there exists $\gam\in\Lam_*(\bfM)$
such that
\eq{support}
\supp f\subset \bigcup\limits_{\gam'\in\gam+\Lam_*(\bfM)^+} X_P^{\gam'}.
\end{equation}
\conj{boundedsupport}
Let $f\in\calS(G,M)$. Then for every $P\in \calP(M)$
the function $f_P$ on $X_P$ has bounded support.
\econj
We don't know how to prove this conjecture in general. In \cite{BK99} we proved it for
$M=T$. It follows also from \reft{mainkinvariant} below that \refco{boundedsupport} holds also
for $K$-invariant elements in $\calS(G,M)$.
\ssec{}{$\calS(G,M)$ as a module over $\calH(M^\ab)$}
Let $\calH(M^\ab)$ denote the Hecke algebra of $M^\ab$.
We have the natural decomposition
\eq{}
\calH(M^\ab)=\bigoplus\limits_{\sig\in\Hom(\bfM^{\ab}(\calO)\to\CC^\x)}\calH(M^\ab)_\sig.
\end{equation}
For each $\sig$ the algebra $\calH(M^\ab)_\sig$ is non-canonically
isomorphic to the algebra of Laurent polynomials $\CC[t_1,...,t_l,t_1^{-1},...,t_l^{-1}]$ where
$l=\dim\bfM^\ab$.

The space $\calS(G,M)$ is naturally a module over $M^\ab$ and hence over $\calH(M^\ab)$.
Let $\calS(G,M)_\sig$ denote the part of $\calS(G,M)$ on which $\bfM^\ab(\calO)$ acts by means of
the character $\sig$.

The following lemma is never used in the sequel but we think that it gives some intuition about the space $\calS(G,M)$:
\lem{}
$\calS(G,M)_\sig$ is a free module over $\calH(M^\ab)_\sig$
(as a $G$-module).
\elem
\prf
First of all by a theorem of Quillen it is enough to prove that $\calS(G,M)_\sig$ is locally free over
$\calH(M^\ab)_\sig$.

It is clear that for each $P$ the $G$-module $\calC_c(X_P)_\sig$ is free over
$\calH(M^\ab)_\sig$. Let now $\chi:M^\ab\to\CC^\x$ be a character
such that $\chi|_{\bfM^\ab(\calO)}=\sig$. We denote by
$\calS(G,M)_\chi$ (resp. $\calC_c(X_P)_\chi$) the $G$-module
of $(M^\ab,\chi)$-coinvariants on $\calS(G,M)$ (resp. on $\calC_c(X_P)$).
It is now easy to see that for every $\chi$ there exists
$P\in\calP(M)$ such that the inclusion
$\calC_c(X_P)_{\chi'}\hookrightarrow\calS(G,M)_{\chi'}$
for every $\chi'$ in some neighbourhood of $\chi$ which implies
what we need.
\epr

\ssec{}{$K$-invariant vectors}For every $\gam\in\Lam_*(\bfM)$ we denote by
$\del_{P,\gam}$ the following function on $X_P$:
$$
\del_{P,\gam}(x)=
\begin{cases}
q^{\la \gam,\rho-\rho_M\ra},\quad \text{if $x\in X_P^\gam$}\\
0, \qquad\text{otherwise.}
\end{cases}
$$
Let $L=\oplus L_i$ be a graded representation of $\bfZ(\bfM^{\vee})$.
Assume that for each $\gam\in\Lam_*(\bfM)$ the multiplicity
of $\gam$ in $L$ is finite. Then we define a function
$\phi^L_P\in\calC(X_P)^K$ by setting
\eq{}
\phi^L_P|_{X_P^\gam}=\sum \dim L_i^\gam q^{-i+\la\gam,\rho-\rho_\bfM\ra}
\end{equation}
Similarly, for every $\mu\in\Lam_*(\bfM)$ we define
\eq{}
\phi^L_{P,\mu}|_{X_P^\gam}=\sum \dim L_i^{\gam-\mu} q^{-i+\la\gam,\rho-\rho_\bfM\ra}
\end{equation}
In other words,
$$
\phi^L_{P,\mu}=\sum \dim L_i^{\gam-\mu} q^{-i}\del_{P,\gam}.
$$

Let $e,f,h\in\grm^{\vee}$ be a principal $sl_2$-triple.
Take
$$
L=\text{Sym}((\gru_\grp^{\vee})^e).
$$
Clearly $L$ has a natural action of
$\bfZ(\bfM^{\vee})$. Also the $L$ carries a natural action of
$h$.

Define a grading on $L$ in the following way.
Let $x\in \text{Sym}^k((\gru_\grp^{\vee})^f)$
Assume that $x$ has eigenvalue $j$ with respect to $h$.
Then we say that $x$ has grading $k+j$.

We set $c_{P,\mu}=\phi^L_{P,\mu}$ for $L$ as above.

We would like to understand the structure of
$\calS(G,M)^K$. Since every $K$-orbit in $X_P$ is
$\bfM^\ab(\calO)$-invariant it follows that every
element of $\calS(G,M)^K$ is automatically
$\bfM^\ab(\calO)$-invariant. Thus $M^\ab$-action on
$\calS(G,M)$ reduces to an $M^\ab/\bfM^\ab(\calO)=\Lam_*(\bfM)$-action
on $\calS(G,M)$.
\th{mainkinvariant}
\begin{enumerate}
\item
For every $P$ and $\mu$ as above
we have $c_{P,\mu}\in\calS(G,M)$.
\item
For every $P\in\calP(M)$ the functions
$c_{P,\mu}$ $($where $\mu$ runs over all elements of $\Lam_*(\bfM)$\,$)$
form a basis in $\calS(G,M)^K$.
\item
For every $P,Q\in \calP(M)$ we have
\eq{}
\calF_{P,Q,\psi}(c_{P,\mu})=c_{Q,\mu}.
\end{equation}
\end{enumerate}
\eth
\prf
We argue along the lines of the proof of Theorem 3.13 in \cite{BK99}.

First of all we claim that $c_{P,\mu}\in L^2(X_P)$. Recall that we denote by
$L$ the space $\text{Sym}((\gru_\grp^\vee)^e)$ with the grading discussed above.
For every $\gam\in\Lam_*(\bfM)$
let
$$
\calK_P(\gam)=\sum \dim L_i^\gam  q^{-i}.
$$
Then to prove that $c_{P,\mu}\in L^2(X_P)$ we must show that
the series
\eq{ltwonorm}
\sum_{\gam\in\Lam_*(\bfM)} \calK_P(\gam)^2
\end{equation}
is convergent. However, it is easy to see that the series
$$
\sum_{\gam\in\Lam_*(\bfM)} \calK_P(\gam)
$$
is convergent to $\prod\limits_{i=1}^\infty (1-q^{-i})^{\dim(\gru_\grp^\vee)^e_i}$.
Hence \refe{ltwonorm} is convergent too.

Since $c_{P,\mu}\in L^2(X_P)$ it follows that $\calF_{P,Q,\psi}(c_{P,\mu})$
is well-defined. Let us show that point 3 of \reft{mainkinvariant} implies points
1 and 2.

 We have the natural isomorphism
$\CC[\bfZ(\bfM^\vee)]\simeq \CC[\Lam_*(\bfM)]$ where $\CC[\bfZ(\bfM^\vee)]$ denotes the algebra
of regular functions on $\bfZ(\bfM^\vee)$ and $\CC[\Lam_*(\bfM)]$ denotes the group
algebra of $\Lam_*(\bfM)$. Let $\calW$ denote the span of $c_{P,\mu}$
for all $\mu\in\Lam_*(\bfM)$. It follows from 3 that $\calW$ does not depend
on $P$ as a subspace of $L^2(G,M)$. We may identify $\calW$ with $\CC[\bfZ(\bfM^\vee)]$
by identifying $c_{P,\mu}$ with $\mu$ (again, this doesn't depend on $P$).
For every $P\in\calP(M)$ we set $\calV_P=\calC_c(X_P)^K$. To prove 1 and 2
we need to show that $\calV_P\subset \calW$ for every $P$ and that
$\calW=\text{span}\{ \calV_P\}_{P\in\calP(M)}$.

Let $\kap:SL(2)\to \bfM^\vee$ be the homomorphism corresponding to the
$sl_2$-triple $(e,f,h)$ chosen above. Let
$$
H_q=\kap
\begin{pmatrix}
q&0\\
0&q^{-1}
\end{pmatrix}.
$$
Define $\bfd_{\bfP}\in\CC[\bfZ(\bfM^\vee)]$ by
$$
\bfd_\bfP(z)=\det(1-H_q^{-1}z)|_{(\bfu_\bfp^\vee)^e}.
$$
Then by definition we have $\bfd_\bfP(c_{P,\mu})=\del_{P,\mu}$.
Thus $\del_{P,\mu}\in \calW$, hence $\calV_P\subset \calW$ for
every $P$. Moreover, as a subspace of $\CC[\bfZ]$ the space
$\calV_P$ is equal to the ideal generated by $\bfd_\bfP$.
Applying Hilbert Nullstellensatz we see that points 1 and 2  of
\reft{mainkinvariant}
follow from the following lemma whose proof is left to the reader.
\lem{}
For every $z\in \bfZ(\bfM^\vee)$ there exists $P\in\calP(M)$ such that
$\bfd_\bfP(z)\neq 0$.
\elem
Let us now prove 3. For this it is enough to show that
\eq{delta}
\calF_{P,Q,\psi}(\del_{P,\mu})=\frac{\bfd_\bfQ}{\bfd_\bfP}\del_{Q,\mu}.
\end{equation}
Define
$$
{\widetilde\bfd_\bfP}(z)=\det(1-H_q^{-1}z)|_{\bfu_\bfp^\vee}\quad\text{and}\quad
\bfd_\bfP^{\text{junk}}=\det(1-H_q^{-1}z)|_{\gru_\grp^\vee/(\gru_\grp^\vee)^e}=
\frac{{\widetilde\bfd_\bfP}(z)}{\bfd_\bfP(z)}.
$$
It follows from fromula 3.22 in \cite{BK99} that
$$
\calG_{P,Q,\psi}(\del_{P,\mu})=\frac{\widetilde\bfd_\bfQ}{\widetilde\bfd_\bfP}\del_{Q,\mu}.
$$
Thus to prove \refe{delta} we need to show that
$\eta_{P,Q,\psi}^{\text{junk}}(\del_{P,\mu})=\frac{\bfd_\bfP^{\text{junk}}}{\bfd_\bfQ^{\text{junk}}}\del_{P,\mu}$
which is left  to the reader.
\epr

In particular we see that for every $P,Q\in\calP(M)$ we have
$$
\calF_{P,Q,\psi}(c_{P,0})=c_{Q,0}.
$$
In other words there exists canonical element $c\in \calS(G,M)$ such that
$c_P=c_{P,0}$ for every $P$.

\bigskip
\noindent
{\it Remarks.} 1. It follows from Theorem 7.3 of
\cite{BFGM} that the function $c_{P,0}$ is equal
(in the appropriate sense) to the function obtained by
traces of Frobenius on the stalks of the intersection
cohomology sheaves of the Drinfeld compactification
$\overline{\text{Bun}_P}$ of the moduli stack
of $P$-bundles on a smooth projective algebraic curve
over $\FF_q$. This is certainly not unexpected.
We don't know, however, how to prove this directly.
Note, however that the spaces $\overline{\text{Bun}_P}$
(or at least their singularities) have their local
counterparts because of \cite{GrKa} and \cite{Dr}.

2. In \cite{BK99} we also study the space $\calS(G,T)^I$
where $I\subset G$ is an Iwahori subgroup.
We give two interpretations of this space: one using the
{\it periodic Hecke module} introduced in \cite{Lu80} and the other
using certain equivariant $K$-group. Similar descriptions
are probably possible for the space $\calS(G,M)^I$ (the corresponding
parabolic periodic Hecke module is studied in \cite{Lu97} and its
interpretation using $K$-theory is given in \cite{LuK1} and \cite{LuK2}).
We have not checked this precisely.

\sec{global}{The case of global fields}
\ssec{}{Remarks on archimedian places}
Let now $\bfG$ be as before but assume that the field $F$ is archimedian.
The operators $\calF_{P,Q,\psi}$ are defined in this case in the same
way as for non-archimedian $F$ (it is not difficult to adjust the proof of unitarity
of $\calF_{P,Q,\psi}$ to the archimedian case; we leave it to the reader).

Let $M$ be a Levi subgroup of $P$.
For every $P\in\calP(M)$ we define $\calS(X_P)$ in the following
way. Let $\bfX_\bfP^\RR$ denote the algebraic variety over $\RR$
obtained by restriction of scalars from $F$ to $\RR$. Let $\calD(\bfX_\bfP^\RR)$
denote the algebra of differential operators on $\bfX_\bfP^\RR$ with polynomial
coefficients. Clearly $\calD(\bfX_\bfP^\RR)$ acts on $C^\infty(X_P)=C^\infty(\bfX_\bfP(F))=
\calD(\bfX_\bfP^\RR(\RR))$. We define
$$
\calS(X_P)=\{ f\in C^\infty(X_P)|\ d(f)\in L^2(X_P)\quad
\text{for every $d\in\calD(\bfX_\bfP^\RR)$}\}.
$$
\conj{real}
For every $P,Q\in\calP(M)$ we have
$$
\calF_{P,Q,\psi}(\calS(X_P))=\calS_(X_Q).
$$
\econj

\noindent
{\bf Example.} Let $G=\Sl(n,\RR)$ and let $P, Q$ be as in \refss{intro-example}. In this case
we have $\calS(X_P)=\calS(V)=$ the Schwartz space of rapidly decreasing functions on $V$. Also
$\calS(X_Q)$ is the Schwartz space of $V^*$. Thus in this case the above conjecture reduces to the
well-known fact saying that
the Fourier transform maps $\calS(V)$ to $\calS(V^*)$.

Assuming that the conjecture is true we may define $\calS(G,M)=\calS(X_P)$ for any $P\in\calP(M)$
(\refco{real} thus says that $\calS(G,M)$ does not depend on the choice of $P$).

We do not know how to prove this conjecture in general. In the sequel we shall choose any
$G\x M^\ab$-invariant subspace
of $L^2(G,M)$ containing $C^\infty_c(X_P)$ for every $P$ and consisting of $G\x M^\ab$-smooth
vectors. We shall denote this subspace by $\calS(G,M)$.

\ssec{}{The global space $\calS_K(G,M)$}
Let $K$ be a global field and let $\calV(K)$ denote its set of places. We now assume that we are given
a reductive group $\bfG$ with a Levi subgroup $\bfM$ as above, both defined over $K$.
For every $v\in\calV(K)$ we
let $K_v$ denote the corresponding completion of $K$. We define
$\calS_v(G,M)=\calS(\bfG(K_v),\bfM(K_v))$. For every finite place $v$ we have canonical
spherical function $c_v\in\calS_v(G,M)$.

Let now $\calS_K(G,M)$ denote the restricted tensor product of the spaces $\calS_v(G,M)$ with
respect to the functions $c_v$.

\th{poisson}There exists a unique $\bfG(K)\x \bfM^\ab(K)$-invariant functional $\eps$ on the space
 $\calS_K(G,M)$
such that the following condition holds:

Let $f\in \calS_K(G,M)$ and assume that for some $v\in \calV(K)$ and $\bfP\in\calP(\bfM)$
the function $f$ lies in the tensor product
$$
\calC_c(\bfX_\bfP(K_v))\ten \bigotimes\limits_{v'\neq v} \calS_{v'}(G,M).
$$
Then
\eq{poisson}
\eps(f)=\sum\limits_{x\in \bfX_\bfP(K)} f_P(x)
\end{equation}
\eth

\prf
Arguing as in Section 6 of \cite{BK99} we see that it is enough to prove the following statement:
let $v_1,v_2$ be two places of $K$ and assume that
$$
f\in \calC_c(\bfX_\bfP(K_{v_1}))\ten \bigotimes\limits_{v'\neq v_1} \calS_{v'}(G,M)\cap
\calC_c(\bfX_\bfQ(K_{v_2}))\ten \bigotimes\limits_{v'\neq v_2} \calS_{v'}(G,M).
$$
Then
$$
\sum\limits_{x\in \bfX_\bfP(K)} f_P(x)=\sum\limits_{x\in \bfX_\bfQ(K)} f_Q(x).
$$
By using exactly the same argument as in Section 7.22 of \cite{BK99} we may see that this is equivalent to the functional
equation for Eisenstein series (induced respectively from characters of $P$ and $Q$).
\epr
\sec{L-function}{Connection with $L$-functions for classical groups}
In this section we indicate how one can use the above results in order to give construction
of (and prove the standard properties) of L- and $\eps$-functions associated with
the standard representation of a classical group. The details will appear elsewhere.
These L-functions were studied
in the (nowadays classical) work \cite{GPSR}. The advantage of our approach is that it may be viewed
as a ``direct'' generalization of the work of Godement and Jacquet (\cite{GJ}) where
the case of $\Gl(n)$ is studied. A similar approach is discussed in
\cite{PSR}. The main ingredient which makes the presentation of this paper
different from  \cite{GPSR} and \cite{PSR}  is
the space $\calS(G,M)$ which was missing in {\it loc. sit}. Because of this in \cite{GPSR} the local zeta-integrals gave
the $L$-function divided by certain auxilliary denominator (which was equal to some product
of abelian $L$-functions). This denominator is absent in our formulation.
\ssec{}{}
Let $\bfH$ denote one of the groups $\Gl(n)$, $\Sp(2n)\x \GG_m$,
or $\Gspin(n)$ (by the definition $\Gspin(n)$ is the quotient of
$\Spin(n)\x \GG_m$ by the diagonal copy of central $\ZZ_2$). To any $\bfH$
like that we associate another reductive
group $\bfG$ together with a maximal  parabolic subgroup $\bfP$ in the following
way.

1. If $\bfH=\Gl(n)$ we set $\bfG=\Sl(n^2)$, $\bfP$ is the stabilizer of a line in the standard representation.

2. If $\bfH=\Sp(2n)\x \GG_m$ we set $\bfG=\Sp(4n)$ and take $\bfP$  to consist of all matrices that stabilize
a Lagrangian subspace.

3. If $\bfH=\Gspin(n)$ we set $\bfG=\Spin(2n)$ and take $\bfP$ to be the stabilizer of a maximal
isotropic subspace in the standard representation.

In all these cases we denote by $\bfM$ the corresponding Levi subgroup of $\bfG$.

For $\bfH$ as above there is a natural character $\sig:\bfH\to\GG_m$ (in case 1 we have
$\sig=\det$; in cases 2 $\sig$ is just the projection to the second multiple; in case 3 if
$g\in \Gspin(n)$ is the image of an element $(g',t)$ in $\Spin(n)\x \GG_m$ then we
let $\sig(g)=t^2$).

Let ${\widetilde\bfH^2}$
denote the set of all pairs $(h_1,h_2)\in\bfH^2$ such that
$\sig(h_1)=\sig(h_2)$. Let $\bfH_1$ denote the kernel of $\sig$ in $\bfH$.
We have the natural diagonal embedding $\Del:\bfH_1\to{\widetilde\bfH^2}$.
Hence $\bfH_1$ acts naturally on ${\widetilde \bfH^2}$ (on the right).
Moreover the quotient ${\widetilde \bfH^2}/\bfH_1$ is naturally isomorphic
to $\bfH$.

We claim now that in all of the above cases there exists an embedding
$\eta:{\widetilde \bfH^2}\hookrightarrow \bfG\x {\bfM^\ab}$ such that the following
hold. Let
$$
\tbfP=\{ (g,m)\in \bfP\x \bfM^\ab|\ \text{such that the image of $g$ in $M^\ab$ is equal to $m$}\}
$$
Thus we have the natural isomorphism
$$
\bfG/[\bfP,\bfP]\simeq \bfG\x \bfM^\ab/\tbfP.
$$
Then we require that

(i) $\eta^{-1}(\tbfP)=\bfH_1$.

(ii) The resulting map $\thsq/\bfH_1\hookrightarrow \bfG/[\bfP,\bfP]$
is an open embedding.

Let us explain the construction of $\eta$.
From now on we shall present all the constructions in cases 1 and 2
only. Case 3 is always similar to Case 2 and we leave it to the reader.

\medskip

\noindent
{\it Case 1.} In this case we identify $\bfG$ with $\Sl(\bfM_n)$ where
$\bfM_n$ denotes the space of $n\x n$-matrices and we take $\bfP$ to be the
stabilizer of the line spanned by the identity matrix.
Then for all $(g_1,g_2)\in \thsq$ the projection of $\eta(g_1,g_2)$ to
$\Sl(\bfM_n)$ takes every $x\in \bfM_n$ to $g_1xg_2^{-1}$.
The composition of $\eta$ with the projection to $\bfM^\ab\simeq \GG_m$
is equal to $\deg g_1^{-1}=\det g_2^{-1}$. In this case
we have $\bfX_\bfP=\bfM_n\backslash\{ 0\}$ and the resulting embedding
$\bfG\hookrightarrow \bfX_\bfP$ is the natural one.

{\it Case 2.} Let $\bfH=\Sp(W,\ome)\x\GG_m$ where $(W,\ome)$ is a symplectic vector
space. Set $V=W\oplus W$ and equip it with the symplectic
form $(\ome\oplus (-\ome))$.  Let $\Del W$ be the diagonal copy of $W$ in $V$.
This is clearly a Lagrangian subspace of $V$.
Then we may identify $\bfG$ with $\Sp(V)$ and take $P$ to be the stabilizer
of $\Del W$. We define
$$
\eta(g_1,g_2,t)=((g_1,g_2),t)
$$
where $(g_1,g_2)$ is an element of $\Sp(V)$ given by
$$
(g_1,g_2)(w_1,w_2)=(g_1(w_1),g_2(w_2)).
$$

The definition of $\eta$ in Case 3 is analogous to Case 2 and we leave it to
the reader.

\ssec{}{The space $\calS(H)$}
We define $\calS(H)=\calS(G,M)$. It follows from \refl{locallyconstant}
that every $f\in \calS(H)$ is locally constant.

Let us consider Case 1 above, i.e. the case $\bfH=\Gl(n)$. In this case
\refss{glnexample} implies that $\calS(H)$ is equal to the space of
locally constant compactly supported functions on $M_n$.
\ssec{}{The ``Fourier transform'' $\calF_{H,\psi}$}
Let $\obfP$ be a parabolic subgroup opposite to $\bfP$.
We claim that in all the cases 1,2,3 above we can identify
$\bfX_\bfP$ with $\bfX_{\obfP}$. Indeed, in Case 1
we just need to identify $\bfM_n$ with the dual vector space.
This can be done by means of the standard bilinear form
$$
(A,B)= \tr AB.
$$
In Case 2 both $\bfX_\bfP$ and $\bfX_{\obfP}$ can be identified with the variety of
of Lagrangian subspaces of $\bfV$ equipped with a volume form (note that in this case
$\bfP$ and $\obfP$ are conjugate to each other). Case 3 is treated similarly.

It is easy to see that  this choice is compatible with the invariant measures on
$X_P$ and $X_{\oP}$. Hence we may view $\calF_{P,\oP,\psi}$ as an operator
acting from $L^2(X_P)$ to $L^2(X_P)$.

\noindent
{\bf Example}. Consider Case 1. Then $\calF_{P,\oP,\psi}$ is just the Fourier transform operator
on $M_n$. Namely given a function $f\in L^2(M_n)$ we have
$$
\calF_{P,\oP,\psi}(f)(x)=\int\limits_{M_n} f(y)\psi(\tr(xy))dy.
$$
By the definition we have
$$
\calS(H)=\calC_c(X_P)+\calF_{P,\oP,\psi}(\calC_c(X_{\oP}))\subset L^2(X_P).
$$
Hence $\calF_{P,\oP,\psi}$ acts from $\calS(H)$ to $\calS(H)$. We denote this operator
by $\calF_{H,\psi}$.
The proof of the following lemma is left to the reader.
\lem{}
With the above identifications we have
$\calF_{\oP,P,\psi}=\calF_{H,\psi^{-1}}$.
\elem
The lemma implies that $\calF_{H,\psi}$ and $\calF_{H,\psi^{-1}}$ are inverse to each other.
Let now $\pi$ be an irreducible representation of $H$. The following theorem
is proved in \cite{GJ} in Case 1 above. In other cases the proof may be obtained
by a similar analysis. We shall not write the details here.
\th{L-function}
\begin{enumerate}
\item
Let $m$ be any matrix coefficient of $\pi$ and let $f\in \calS(H)$
be any function.
Then the integral
$$
Z(f,m,s)=\int\limits_H m(h)f(h)|\sig(h)|^s dh
$$
is absolutely convergent for $Re(s)>>0$.
\item
$Z(f,m,s)$ extends meromorphically to the whole of $\CC$ as a rational
function of $q^{-s}$.
\item
For a fixed $\pi$ the functions $Z(f,m,s)$ form a fractional ideal $J_\pi$
in the ring $\CC[q^{-s}]$ of polynomials $q^{-s}$.
\item
There exists unique polynomial $P_\pi(t)$ with $P_\pi(0)=1$ such that
$J_\pi$ is generated by $\frac{1}{P_\pi(q^{-s})}$. We set
$$
L(\pi,s):=\frac{1}{P_\pi(q^{-s})}.
$$
\item
In Cases 2 and 3
$L(\pi,s)$ coincides with the local $L$-function of $\pi$ constructed in
\cite{GPSR} (in Case 1 the above definition of $L(\pi,s)$ is exactly the same
as the one in \cite{GJ}).
\end{enumerate}
\eth

One can use the operator $\calF_{H,\psi}$ in order to define the corresponding $\eps$-factors
(using the idea in \cite{GJ} where this is done in Case 1).
Also, one can show that the generalized Poisson summation formula \refe{poisson}
implies that the global version of $L(\pi,s)$ satisfies the corresponding
functional equation.

\end{document}